\documentclass{article}

\usepackage{PRIMEarxiv}
\usepackage[utf8]{inputenc} 
\usepackage[T1]{fontenc}    
\usepackage{hyperref}       
\usepackage{url}            
\usepackage{booktabs}       
\usepackage{amsfonts}       
\usepackage{nicefrac}       
\usepackage{microtype}      
\usepackage{lipsum}
\usepackage{fancyhdr}       
\usepackage{graphicx}       
\graphicspath{{media/}}     
\usepackage{amsmath,amssymb, amsthm}

\title{New inverse problems for a time-switched system of wave and diffusion equations
}

\author{
  Karimov E.T.\footnote{Corresponding author} \\
  Department of Mathematics: Analysis, Logics and Discrete Mathematics\\
  Ghent University\\ 
  Ghent \\
    \texttt{erkinjon.karimov@ugent.be} \\
   \And
  Murolimova N.A. \\
  Department of Mathematical Analysis and Differential Equations \\
  Fergana State University \\
  Fergana\\
  \texttt{nmurolimova@gmail.com} \\
}

\begin{document}

\maketitle

\begin{abstract}
We study two new classes of inverse problems for a time-switched system in which a
fractional wave equation (with Caputo derivative of order $\alpha \in (1,2)$) governs the
dynamics on the interval $[0,a)$, and a fractional diffusion equation (with Caputo derivative
of order $\beta \in (0,1)$ taken with respect to the switching point $t=a$) governs the
dynamics on $(a,b]$.
The two problems differ in which part of the transmitting condition at the interface $t=a$ is
regarded as unknown. In both cases the overdetermination data consist of a single spatial measurement of the
solution at a fixed time $\xi \in (a,b)$.
Using the spectral expansion method with respect to the classical Sturm-Liouville
eigensystem on $[0,1]$, we reduce each problem to a sequence of coupled scalar Cauchy
problems involving the two-parameter Mittag-Leffler function.
Explicit series representations for the solution $u(t,x)$ and the unknown interface
functions $h(x)$ and $\bar{h}(x)$ are derived.
Uniform convergence of the resulting infinite series and their relevant derivatives is
established through four auxiliary lemmas, using the decay estimates for the Mittag-Leffler
function, integration-by-parts arguments, the Cauchy--Schwarz inequality, and the Weierstrass
$M$-test.
A uniqueness and existence theorem is stated for Problem~1 under explicit Sobolev-type
regularity conditions on the data, with an analogous result for Problem~2.
\end{abstract}
 
 
 
\section{Introduction}
 
\subsection*{Background and motivation}
 
Fractional-order partial differential equations have become indispensable tools for modelling
physical processes that exhibit memory effects, anomalous transport, or intermediate behaviour
between classical diffusion and wave propagation.
The foundational well-posedness and spectral theory for initial-boundary value problems
governed by the Caputo time-fractional diffusion-wave equation
\[
  \partial_t^\alpha u - Lu = f, \qquad 0 < \alpha \leq 2,
\]
where $L$ is a uniformly elliptic spatial operator, was rigorously established by
Sakamoto and Yamamoto~\cite{SakamotoYamamoto2011}.
Their eigenfunction-expansion framework, together with decay estimates for the Mittag-Leffler
function has since become a standard starting point for the analysis of both direct and
inverse problems in this setting; we follow the same approach in the present work.
 
A distinctive and practically important class of models arises when the governing equation
is not fixed throughout the time interval, but \emph{switches} between two different
fractional orders at a prescribed transition time $t = a$.
Such \emph{time-switched} or \emph{composite-type} systems naturally appear whenever a
physical or biological process undergoes a qualitative change of regime at a critical instant.
Two canonical scenarios motivate the present paper.
 
\paragraph{Scenario 1: Thermoelastic and viscoelastic transitions.}
In many solid materials, mechanical behaviour transitions from wave-like (hyperbolic,
$1 < \alpha < 2$) at short time scales to diffusion-like (parabolic, $0 < \beta < 1$) at
longer time scales as internal energy redistribution takes over from elastic restoring forces.
Models of anomalous thermoelastic relaxation, polymer viscoelasticity, and
fractional Maxwell fluids exhibit precisely this structure (see, e.g.,
\cite{MetzlerKlafter2000,Mainardi2010}).
The transmitting conditions at $t = a$ encode the physical balance at the
moment of transition: continuity of displacement (the position condition) and continuity
or prescribed jump of the velocity-type flux (the velocity condition).
In practice, the flux exchange at the transition instant may depend on the internal state of
the medium and cannot always be observed directly, rendering its spatial distribution an
unknown to be recovered -- exactly the problem treated here.
 
\paragraph{Scenario 2: Phase-change diffusion and anomalous transport in porous media.}
Contaminant transport and thermal energy redistribution through heterogeneous porous media
often exhibit a two-stage behaviour: an initial inertia-dominated phase (superdiffusion,
captured by a fractional order $\alpha > 1$) followed by a subdiffusive regime (fractional
order $\beta < 1$) once the medium becomes locally saturated or once the driving pressure
gradient subsides.
The interface between these two regimes is a temporal interface analogous to a Stefan-type
free boundary, but now occurring in the time variable rather than in space.
At the switching instant, the concentration profile satisfies a continuity condition, while the flux (or its fractional analogue) may experience a spatially distributed jump due to
adsorption, trapping, or source activation at $t = a$.
If this jump is not measured directly -- a situation common in subsurface monitoring -- its
recovery from a single spatial snapshot at a later time $\xi > a$ constitutes the inverse
problem studied in Problem~1.
Conversely, if the flux is prescribed but the concentration undergoes a sudden spatially
non-uniform offset (e.g., a localised source activation or a chemical transformation at
$t = a$), then Problem~2 is the relevant formulation.
 
\paragraph{Summary of the novelty.}
The identification of unknown \emph{interface conditions} -- as opposed to unknown source
terms, coefficients, or initial data -- constitutes a genuinely new class of inverse problems
for time-fractional evolution equations.
To the best of our knowledge, this formulation has not been considered previously in the
fractional calculus literature.
It is to be distinguished from the extensive existing theory of spatial-interface (transmission)
inverse problems for parabolic and hyperbolic PDEs, where the interface is a \emph{spatial}
surface separating two media with different coefficients (see, e.g.,
\cite{BellassoueYamamoto2006,WenZhang2022}).
In the present work the interface is \emph{temporal}, and the unknowns $h(x)$ and $\bar{h}(x)$ describe the mismatch in the dynamic transmitting conditions at the switching instant $t = a$
rather than a discontinuity in the spatial medium.
 
\subsection*{Related work}
 
\paragraph{Direct problems for time-fractional equations.}
Since the work of Sakamoto and Yamamoto~\cite{SakamotoYamamoto2011}, the initial-boundary
value theory for fractional diffusion-wave equations has been developed in many directions,
including multi-term operators~\cite{LiLiuYamamoto2015}, time-dependent
coefficients~\cite{KubicaYamamoto2018}, equations on metric
graphs~\cite{KarimovHilfer2024MMAS}, non-self-adjoint operators \cite{BKM26}, and abstract Hilbert-space
settings~\cite{KubicaRyszewskaYamamoto2020,DiffWaveAbstract2025}.

Direct problems for equations combining wave and diffusion dynamics have also been extensively studied. For example, Al‑Salti, Karimov, and Al‑Ghabsi \cite{AKG25} analyzed an initial-boundary value problem for a mixed fractional wave equation, establishing well-posedness and qualitative properties of solutions. Moreover, Karimov and Toshpulatov \cite{KTG25} considered a mixed wave–diffusion–wave equation and proved solvability results for the corresponding initial-boundary value problem. These works provide a solid analytical foundation for the study of inverse problems in hybrid models exhibiting regime switching.
 
\paragraph{Inverse source and coefficient problems.}
Inverse problems for single-regime fractional diffusion-wave equations have been investigated
extensively.
Sakamoto and Yamamoto~\cite{SakamotoYamamoto2011mcrf} studied the recovery of a spatially
varying source by final-time overdetermination.
Cheng and Li~\cite{ChengLi2021} established uniqueness and Lipschitz stability for
spatially dependent sources in non-symmetric fractional diffusion-wave equations.
Durdiev and co-authors~\cite{DurdievRahmonov2023} treated inverse coefficient problems for
the fractional wave equation with Hilfer derivative.
For a systematic survey, we refer to the review~\cite{LiLiuYamamoto2019review}.
 
\paragraph{Inverse problems for mixed-type and multi-term fractional equations.}
The present paper belongs to a line of research initiated by one of the authors on inverse
problems for equations of mixed type involving Caputo fractional derivatives.
Karimov and co-authors established well-posedness results for inverse source problems in
mixed-type integral-differential equations~\cite{KTU} and for
direct and inverse problems in equations combining two different Caputo orders on a spatial
star graph~\cite{KarimovHilfer2024MMAS}.
The extended abstract~\cite{KarimovTokenbetovKhasanov2024} proposed a program for inverse
problems in time-fractional mixed equations, which the present work carries forward.
Inverse source problems for degenerate time-fractional PDEs, where both elliptic degeneracy
and fractional time memory appear simultaneously, were studied by Karimov and
Al-Salti~\cite{KarimovAlSalti2019}.
Salakhitdinov and Karimov~\cite{SalakhitdinovKarimov2016} proved uniqueness of an inverse
source problem for a mixed-type equation via spectral methods.

In addition to inverse source and coefficient identification problems, recent works have addressed inverse problems involving the determination of unknown temporal parameters in mixed wave–diffusion models. In particular, Karimov and Mirzaeva \cite{KM24} investigated an inverse problem for determining an unknown time moment in a time-fractional mixed wave–diffusion equation, establishing uniqueness and solvability results via spectral methods. Furthermore, Karimov and Tokmagambetov \cite{KT24} considered a related inverse problem for identifying an unknown time-switching instant in a mixed wave–diffusion equation, providing rigorous analysis of existence and uniqueness. These studies demonstrate the growing interest in inverse problems where the transition time between regimes itself becomes an unknown parameter, complementing the present work where unknown interface conditions are recovered.

\paragraph{Transmission (spatial-interface) inverse problems.}
In the classical integer-order setting, Bellassoued and
Yamamoto~\cite{BellassoueYamamoto2006} established Lipschitz stability for an inverse source
problem in a parabolic transmission problem with a spatial interface, using Carleman estimates
adapted to the discontinuous setting.
Wen and Zhang~\cite{WenZhang2022} treated the simultaneous recovery of initial data and
source strength in a similar parabolic transmission problem, combining Carleman estimates
with logarithmic convexity and a numerical iterative thresholding scheme.

The present problems differ fundamentally from these works: here the interface is temporal,
the governing operator on each side carries a Caputo fractional derivative of different
order, and the unknowns are components of the transmitting condition itself, not sources
or initial data.
 
\paragraph{Anomalous diffusion, wave--diffusion crossover, and time-switching.}
From the modeling perspective, wave-to-diffusion crossover in time is linked to
continuous-time random walk (CTRW) theory and to the theory of anomalous transport
(see~\cite{MetzlerKlafter2000}).
The mathematical framework of the Heaviside-switched operator employed here,
\[
  H(a-t)\,{}^C\!D_{0t}^\alpha u + H(t-a)\,{}^C\!D_{at}^\beta u - u_{xx} = f,
\]
is a natural fractional analogue of composite-type (parabolic-hyperbolic) equations that have
been studied classically, and it allows the precise characterization of the regularity class
required for solutions to exist and for the inverse problem to be well-posed.
 
\subsection*{Organization of the paper}
 
The remainder of the paper is structured as follows.
In Section~1 we formulate and solve Problem~1 (unknown velocity-type interface function
$h(x)$), derive the explicit series representation of the solution and the unknown, and
prove the uniform convergence of the relevant infinite series through Lemmas~1-4 and
Theorem~1.
Section~2 is devoted to Problem~2 (unknown position-type interface function $\bar{h}(x)$),
which is analyzed by an analogous spectral method; the corresponding Fourier coefficients
$B_k$ and $\bar{h}_k$ are found in closed form.
 
  Throughout the paper, $H(\cdot)$ denotes the Heaviside step function,
  ${}^C\!D_{0t}^\alpha$ and ${}^C\!D_{at}^\beta$ are the Caputo fractional derivatives of
  orders $\alpha$ and $\beta$ taken from the left endpoints $0$ and $a$ respectively, and
  \[
    E_{a,b}(z) = \sum_{n=0}^{\infty} \frac{z^n}{\Gamma(an+b)}, \qquad a,b,z \in \mathbb{C},
    \quad \operatorname{Re}(a) > 0,
  \]
  denotes the two-parameter Mittag-Leffler function.

\section{Inverse problem on finding the unknown part of the velocity transmitting condition.}

\textbf{Problem 1.} We aim to find functions $u(t,x) \in C(\bar{\Omega})\cap C^2(\Omega_1\cup\Omega_2)$, $u(\cdot,x) \in AC^2[0,a)\cup AC^1[a,b),$\, $h(x)\in C[0,1]$, which satisfies the equation
\begin{equation}\label{1}
   H(a-t)\,	{}^CD_{0t}^{\alpha} u(t,x) +	H(t-a)\, ^CD_{at}^\beta u(t,x)-u_{xx}(t,x)=f(t,x) 
\end{equation}
together with the initial condition
\begin{equation}\label{2}
   u(0,x)=\varphi(x), \quad 0\le x\le 1, 
\end{equation}
boundary conditions
\begin{equation}\label{3}
 u(t,0)=0,\,\,u(t,1)=0,\,\,\,0\le t\le b,   
\end{equation}
the over-determination condition for fixed $\xi\in (a,b)$
\begin{equation}\label{4}
 u(\xi,x)=\psi(x), \quad 0\le x\le 1,   
\end{equation}
and the transmitting condition
\begin{equation}\label{5}
    \lim_{t \rightarrow a+0} { }^C D_{at}^{\beta} u(t, x)=u_t(a-0,x)+h(x), \quad 0< x< 1.
\end{equation}
Another transmitting condition follows from the regularity condition $u(t,x) \in C(\bar{\Omega})$, i.e.
\begin{equation}\label{6}
    u(a+0,x)=u(a-0,x), \quad 0\le x\le 1.
\end{equation}
Here $H(\cdot)$ is a Heaviside function, $ \Omega=\Omega_1\cup \Omega_2\cup \{t=a,\,0\le x\le 1\},$ $\Omega_1= \left\{(t,x): 0<x<1,  0<t<a\right\},$ $\Omega_2= \left\{(t,x): 0<x<1,  a<t<b \right\}$, $\varphi(x), \psi(x), f(t,x)$ are given functions and with $a, b, \alpha,  \beta \in \mathbb{R}$ such that $b>a>0, \, 0<\beta<1, \, 1<\alpha < 2$.

To find the series form of the solution to the formulated problem, we use a method of spectral expansion. The spectral problem associated with the \eqref{1} and boundary conditions \eqref{3} is the following:
\begin{equation}\label{7}
    X^{\prime \prime}(x)+\lambda X(x)=0, \quad X(0)=0, \quad X(1)=0 .
\end{equation}

 It is well-known that eigenvalues of this Shturm-Liouville problem \eqref{7} are $\lambda_k=(k \pi)^2$ and corresponding eigenfunctions are $X_k=\sin k \pi x$, and they form a complete orthonormal system in $L_2(0,1)$. Therefore, we search for a solution to problem 1 as follows:
 \begin{equation}\label{8}
     u(t, x)=\sum_{k=1}^{\infty} u_k(t) \sin k \pi x,
 \end{equation}
 where $u_k(t)$ are unknown functions to be found.  Representing an unknown function $h(x)$ as
\begin{equation} \label{h_k}
	h(x)=\sum \limits_{k=1}^\infty h_k \sin k\pi x, \quad 0\le x\le 1,
\end{equation}
we also expand the given function $f(t, x)$ to the same system $\{X_k(x)=\sin k\pi x\}_{k=1}^\infty$.
 \begin{equation}\label{9}
     f(t, x)=\sum\limits_{k=1}^{\infty} f_k(t) \sin k\pi x,
 \end{equation}
 $f_k(t)$ are the Fourier coefficients, defined as 
 $$
 f_k(t)=2 \int_0^1 f(t, x)\sin k\pi x dx.
 $$ 
 We substitute \eqref{8} and \eqref{9} into \eqref{1} and will get
 \begin{equation}\label{10}
     { }^C D_{0t}^\alpha u_k (t)+\lambda_k u_k(t)=f_k(t),\,\,  1<\alpha<2,\,\,  0<t<a
 \end{equation}
 \begin{equation}\label{11}
     { }^C D_{a t}^\beta u_k (t)+\lambda_k u_k(t)=f_k(t),\,\, 0<\beta<1,\,\, a<t<b.
 \end{equation}
 
 The initial condition \eqref{2} gives us the following
 \begin{equation}\label{12}
     u_k(0)=\varphi_k,
 \end{equation}
 where 
 $$
 \varphi_k=2 \int_0^1 \varphi(x) \sin k\pi x d x
 $$ 
 are Fourier coefficients of the function $\varphi(x)$ expanded to the above-mentioned sine-Fourier series. 
 
 The solution of the Cauchy problem \eqref{11}, $u_k\left(a\right)=C_k$ on $\left(a, b\right)$ and \eqref{10}, \eqref{12},$u_k^{\prime}(0)=B_k$ on $\left(0, a\right)$ will have the form
 \begin{equation}\label{13}
      u_k(t)=C_k E_{\beta, 1}\left[-\lambda_k\left(t-a\right)^{\beta}\right]
+\int\limits_{a}^t(t-z)^{\beta-1} E_{\beta, \beta}\left[-\lambda_k(t-z)^{\beta}\right] f_k(z) d z,\,\,\, a\le t\le b,
 \end{equation}
 
\begin{equation}\label{14}
    u_k(t)=\varphi_k E_{\alpha, 1}\left[-\lambda_k t^\alpha\right]+B_k t E_{\alpha, 2}\left[-\lambda_k t^{\alpha}\right]
+\int\limits_0^t(t-z)^{\alpha-1} E_{\alpha, \alpha}\left[-\lambda_k(t-z)^{\alpha}\right] f_k(z) d z,\, 0\le t\le a,
      \end{equation}
correspondingly. Here $B_k, C_k$ are subject to find. 

For this aim, we use over-determination and transmitting conditions.
Substituting \eqref{13} into the condition $u_k(\xi)=\psi_k$ and considering that ${ E_{\beta, 1}\left[-\lambda_k(\xi-a)^{\beta}\right]} \not=0$ for any $0<\beta <1$, we will find $C_k$ as follows:
\begin{equation}\label{15}
    C_k= \frac{1}{E_{\beta, 1}\left[-\lambda_k(\xi-a)^{\beta}\right]}\left[\psi_k-\int_a^{\xi}(\xi-z)^{\beta-1} E_{\beta, \beta}\left[-\lambda_k(\xi-z)^{\beta}\right] f_k(z) d z\right].
\end{equation}
From the second transmitting condition \eqref{6} and equalities \eqref{13}, \eqref{14} we will find $B_k$ as
\begin{equation}\label{16}
    B_k= \frac{1}{a E_{\alpha, 2}\left[-\lambda_ka^{\alpha}\right]}\bigg[C_k-\varphi_k E_{\alpha, 1}\left(-\lambda_ka^{\alpha}\right)
-\int\limits_0^{a}(a-z)^{\alpha-1} E_{\alpha, \alpha}\left[-\lambda_k(a-z)^{\alpha}\right] f_k(z) d z\bigg].
\end{equation} 

Note that here we considered the fact that for $1<\alpha \leqslant \frac{4}{3}$ the function $E_{\alpha, 2}\left[-\lambda_ka^{\alpha}\right]\neq 0$.

We will use the first transmitting condition \eqref{5}, that is 
$$
\lim\limits_{t\to a+0} { }^C D_{at}^{\beta} u_k(t) = u_k^\prime (a-0)+h_k.
$$
One can easily verify that
\begin{equation*}
\lim\limits_{t\to a+0} { }^C D_{at}^{\beta} u_k(t) = f_k\left(a\right)-\lambda_k u_k\left(a+0\right).
\end{equation*}

Considering $$\frac{d}{dt} E_{\beta, \beta}\left[\lambda t^\beta\right]=\lambda t^{\beta-1} E_{\beta, \beta}\left[\lambda t^\beta\right],\,\,\,\frac{d}{dt}\left(t E_{\beta, 2}\left[\lambda t^\beta\right]\right)=E_{\beta, 1}\left[\lambda t^\beta\right]$$ for $\lambda \in \mathbb{R}, \beta \in \mathbb{R}^+$, one can obtain
\begin{equation*}
\begin{gathered}
u_k^{\prime}(t)=-\lambda_k \varphi_k t^{\alpha-1} E_{\alpha, \alpha}\left[-\lambda_k t^{\alpha}\right] +B_k E_{\alpha, 1}\left[-\lambda_k t^{\alpha}\right]+\int_0^t (t-z)^{\alpha-2} E_{\alpha, \alpha-1}\left[-\lambda_k (t-z)^{\alpha}\right] f_k(z) dz,\,\,a<t<b.
\end{gathered}
\end{equation*}

Considering these equalities, we will find $h_k$ as
\begin{equation}\label{17}
   \begin{gathered}
h_k=f_k\left(a\right)+\lambda_k\varphi_ka^{\alpha-1} E_{\alpha,\alpha}\left[-\lambda_k a^{\alpha}\right]-\lambda_kC_k-B_k  E_{\alpha,1}\left[-\lambda_k a^{\alpha}\right]\\
-\int\limits_0^a (a-z)^{\alpha-2} E_{\alpha, \alpha-1}\left[-\lambda_k (a-z)^{\alpha}\right] f_k(z) dz.
\end{gathered}
\end{equation}

\subsection{Convergence of infinite series}

 Now, considering the well-known estimate
 \begin{equation}\label{MLE}
     |E_{\bar a,\bar b}(z)|\le \frac{c}{1+|z|},\,\,z\le 0, \,\bar b\in\mathbb{R},\,0<\bar a<2
 \end{equation}
 for positive constant $c$ independent of $z$, based on \eqref{15} we get 
        \begin{equation}\label{18}
            |C_k|\le c_1|\psi_k|+c_2\int\limits_a^\xi\frac{(\xi-z)^{\beta-1}|f_k(z)|}{1+\lambda_k(\xi-z)^\beta} d z.
        \end{equation}
Using the estimate \eqref{MLE} and considering \eqref{16}, \eqref{18} we obtain 
\begin{equation}\label{19}
        |B_k|\le c_3|\psi_k|+c_4\frac{|\varphi_k|}{\lambda_k}+c_5\int\limits_0^a \frac{(a-z)^{\alpha-1}|f_k(z)|}{1+\lambda_k(a-z)^\alpha}d z+c_6 \int\limits_a^\xi\frac{(\xi-z)^{\beta-1}|f_k(z)|}{1+\lambda_k(\xi-z)^\beta} d z.    \end{equation}
	    
 Based on \eqref{MLE} and \eqref{18}, from \eqref{13} we deduce
 \begin{equation}\label{u_k1}
|u_k(t)|\le c_7\frac{|\psi_k|}{1+\lambda_k(t-a)^\beta}+\frac{c_8}{1+\lambda_k(t-a)^\beta}\int\limits_a^\xi\frac{(\xi-z)^{\beta-1}|f_k(z)|}{1+\lambda_k(\xi-z)^\beta} d z+c_9 \int\limits_a^t \frac{(t-z)^{\beta-1}|f_k(z)|}{1+\lambda_k(t-z)^\beta}dz,\,a\le t\le b.      
 \end{equation}
 Here $c_i\,\,(i=\overline{1,9})$ are positive constants independent of $t$ and $\lambda_k$.
 
 Finally, considering \eqref{u_k1} based on \eqref{8} we will get
        \begin{eqnarray}\label{20}
            \left| u\left( t,x \right) \right|&\le& \sum\limits_{k=1}^{+\infty } \frac{c_7|\psi_k|}{1+\lambda_k(t-a)^\beta}+\sum\limits_{k=1}^{+\infty } \frac{c_8}{1+\lambda_k(t-a)^\beta}\int\limits_a^\xi\frac{(\xi-z)^{\beta-1}|f_k(z)|}{1+\lambda_k(\xi-z)^\beta} d z\nonumber\\
           & +&\sum\limits_{k=1}^{+\infty } c_9 \int\limits_0^t \frac{(t-z)^{\beta-1}|f_k(z)|}{1+\lambda_k(t-z)^\beta}dz,\,\,(t,x)\in \Omega_2.
        \end{eqnarray}
The estimate of infinite series corresponding to the function $u_{xx}(t,x)$ can be written as follows:
\begin{eqnarray}\label{21}
            \left| u_{xx}\left( t,x \right) \right|&\le& \sum\limits_{k=1}^{+\infty } \frac{c_7\lambda_k|\psi_k|}{1+\lambda_k(t-a)^\beta}+\sum\limits_{k=1}^{+\infty } \frac{c_8\lambda_k}{1+\lambda_k(t-a)^\beta}\int\limits_a^\xi\frac{(\xi-z)^{\beta-1}|f_k(z)|}{1+\lambda_k(\xi-z)^\beta} d z\nonumber\\
           & +&\sum\limits_{k=1}^{+\infty } c_9\lambda_k \int\limits_0^t \frac{(t-z)^{\beta-1}|f_k(z)|}{1+\lambda_k(t-z)^\beta}dz,\,\,(t,x)\in \Omega_2.
        \end{eqnarray}

To prove the uniform convergence of these infinite series, we need the following statements.

Let $f(t,x)$ be a real-valued function defined on $[a,b] \times [0,1]$ with $b > a > 0$. Let $f_k(t)$ denote the spatial Fourier sine coefficients of $f(t,x)$:
$$f_k(t) = 2 \int\limits_0^1 f(t,x) \sin(k\pi x) dx.$$
Assume that $\lambda_k = (k\pi)^2$, $\beta \in (0,1)$, and $\xi \in (a,b)$ is a fixed constant. 
Define the two infinite series:
$$
S_1(t) = \sum_{k=1}^\infty \frac{\lambda_k}{1+\lambda_k(t-a)^\beta} \int\limits_a^\xi \frac{(\xi-z)^{\beta-1}|f_k(z)|}{1+\lambda_k(\xi-z)^\beta} dz,$$
$$
S_2(t) = \sum_{k=1}^\infty \lambda_k \int\limits_a^t \frac{(t-z)^{\beta-1}|f_k(z)|}{1+\lambda_k(t-z)^\beta} dz.
$$

{\bf Lemma 1.} 
If $f(t,x)$ satisfies the boundary conditions $f(t,0) = f(t,1) = 0$ for all $t \in [a,b]$, and $f(t, \cdot)$ is absolutely continuous with its partial derivative $\frac{\partial f}{\partial x} \in L^2(0,1)$, then the series $S_1(t)$ and $S_2(t)$ converge uniformly on the interval $[a,b]$.

{\sc Proof:} We begin by evaluating the integral core present in both series. Consider the generic integral $I_k(y)$ for $y \in (a, b]$:
$$
I_k(y) = \int\limits_a^y \frac{(y-z)^{\beta-1}}{1+\lambda_k(y-z)^\beta} dz.
$$
Using the substitution $u = (y-z)^\beta$, we find $du = -\beta(y-z)^{\beta-1} dz$. The limits of integration map $z=a \to u=(y-a)^\beta$ and $z=y \to u=0$.
$$
I_k(y) = \frac{1}{\beta} \int\limits_0^{(y-a)^\beta} \frac{1}{1+\lambda_k u} du = \frac{1}{\beta \lambda_k} \ln\left(1+\lambda_k(y-a)^\beta\right).
$$
Let $M_k = \sup\limits_{z \in [a,b]} |f_k(z)|$. Taking the absolute value of $S_2(t)$ and applying our evaluated integral yields:
$$
|S_2(t)| \le \sum_{k=1}^\infty \lambda_k M_k I_k(t) = \frac{1}{\beta} \sum_{k=1}^\infty M_k \ln\left(1+\lambda_k(t-a)^\beta\right).
$$
Substituting $\lambda_k = k^2\pi^2$ and knowing that the maximum value on $[a,b]$ occurs at $t=b$:
$$
|S_2(t)| \le \frac{1}{\beta} \sum_{k=1}^\infty M_k \ln\left(1+k^2\pi^2(b-a)^\beta\right).
$$
For large $k$, the term $\ln(1+k^2\pi^2(b-a)^\beta) = O(\ln k)$. Thus, uniform convergence is guaranteed if $\sum\limits_{k=1}^\infty M_k \ln k < \infty$. Applying the same bounding logic to $S_1(t)$, maximizing the external multiplier $\dfrac{1}{1+\lambda_k(t-a)^\beta}$ at $t=a$ (where it equals $1$), yields the same asymptotic condition: $\sum\limits_{k=1}^\infty M_k \ln k < \infty$. 

Now we evaluate $M_k$ using integration by parts. Since $f(t,0) = f(t,1) = 0$:
$$
f_k(t) = \left[ -f(t,x) \frac{\cos(k\pi x)}{k\pi} \right]_0^1 + \frac{1}{k\pi} \int_0^1 \frac{\partial f}{\partial x}(t,x) \cos(k\pi x) dx = \frac{\hat c_k(t)}{k\pi},
$$
where $\hat c_k(t)$ are the Fourier cosine coefficients of $\dfrac{\partial f}{\partial x}$. Let $\hat C_k = \sup\limits_t |\hat c_k(t)|$. Our hypothesis states that $\sum\limits_{k=1}^\infty \hat C_k^2 < \infty$. Therefore, $M_k \le \dfrac{\hat C_k}{k\pi}$.

We now test our required series using the Cauchy-Schwarz inequality:
$$
\sum_{k=1}^\infty M_k \ln k \le \frac{1}{\pi} \sum_{k=1}^\infty \hat C_k \frac{\ln k}{k} \le \frac{1}{\pi} \left( \sum_{k=1}^\infty \hat C_k^2 \right)^{\frac{1}{2}} \left( \sum_{k=1}^\infty \frac{\ln^2 k}{k^2} \right)^{\frac{1}{2}}.
$$
Since $\sum\limits_{k=1}^\infty \dfrac{\ln^2 k}{k^2}$ converges (by the integral test) and $\sum\limits_{k=1}^\infty \hat C_k^2 < \infty$ by hypothesis, the sum converges absolutely. Therefore, $S_1(t)$ and $S_2(t)$ converge uniformly on $[a,b]$. $\blacksquare$

Let $\psi(x)$ be a real-valued function defined on $[0,1]$ and let $\psi_k$ denote its Fourier sine coefficients:
$$
\psi_k = 2 \int\limits_0^1 \psi(x) \sin(k\pi x) dx.
$$
Assume $\lambda_k = (k\pi)^2$ and $\beta \in (0,1)$. Define the infinite series:
$$
S(t) = \sum_{k=1}^\infty \frac{\lambda_k|\psi_k|}{1+\lambda_k(t-a)^\beta}.
$$

{\bf Lemma 2.} 
If $\psi(x)$ lies in the Sobolev space $H^3(0,1)$ - meaning $\psi$, $\psi'$, and $\psi''$ are absolutely continuous and $\psi''' \in L^2(0,1)$ - and satisfies the boundary conditions $\psi(0) = \psi(1) = \psi''(0) = \psi''(1) = 0$, then the series $S(t)$ converges uniformly on $[a,b]$.

{\sc Proof:} By the Weierstrass M-test, the series $S(t)$ converges uniformly on $[a,b]$ if we can bound the terms by a convergent sequence of constants independent of $t$. Because $\lambda_k = k^2\pi^2 > 0$, $\beta > 0$, and $(t-a) \ge 0$ for all $t \in [a,b]$, the denominator satisfies:
$$
1 + \lambda_k(t-a)^\beta \ge 1.$$
This absolute minimum is attained uniquely at $t=a$. Therefore, the supremum of each term in the series over the interval $[a,b]$ is found by setting $t=a$:
$$
\sup_{t \in [a,b]} \left| \frac{\lambda_k|\psi_k|}{1+\lambda_k(t-a)^\beta} \right| = \lambda_k |\psi_k| = k^2\pi^2|\psi_k|.
$$
Thus, uniform convergence is established if the series $\sum\limits_{k=1}^\infty k^2 |\psi_k|$ converges. To bound $\psi_k$, we apply integration by parts to the definition of the Fourier coefficients. Using the condition $\psi(0) = \psi(1) = 0$:
$$
\psi_k = \frac{2}{k\pi} \int\limits_0^1 \psi'(x) \cos(k\pi x) dx.
$$
Applying integration by parts a second time, and noting that the boundary terms evaluate to zero because $\sin(k\pi x) = 0$ at $x=0,1$:
$$
\psi_k = -\frac{2}{k^2\pi^2} \int\limits_0^1 \psi''(x) \sin(k\pi x) dx.
$$
Applying integration by parts a third time, using the given boundary conditions $\psi''(0) = \psi''(1) = 0$:
$$
\psi_k = -\frac{2}{k^3\pi^3} \int\limits_0^1 \psi'''(x) \cos(k\pi x) dx = -\frac{1}{k^3\pi^3} \bar c_k,
$$
where $\bar c_k = 2 \int\limits_0^1 \psi'''(x) \cos(k\pi x) dx$ are the Fourier cosine coefficients of the third derivative $\psi'''(x)$. Since we assumed $\psi''' \in L^2(0,1)$, Bessel's inequality guarantees that $\sum\limits_{k=1}^\infty \bar c_k^2 < \infty$. Substituting our expression for $\psi_k$ into our convergence requirement:
$$
\sum_{k=1}^\infty k^2 |\psi_k| = \sum_{k=1}^\infty k^2 \frac{|\bar c_k|}{k^3\pi^3} = \frac{1}{\pi^3} \sum_{k=1}^\infty \frac{|\bar c_k|}{k}.
$$
Using the Cauchy-Schwarz inequality:
$$
\sum_{k=1}^\infty \frac{|\bar c_k|}{k} \le \left( \sum_{k=1}^\infty \bar c_k^2 \right)^{\frac{1}{2}} \left( \sum_{k=1}^\infty \frac{1}{k^2} \right)^{\frac{1}{2}}.
$$
Because $\sum\limits_{k=1}^\infty \dfrac{1}{k^2} = \dfrac{\pi^2}{6}$ and $\sum\limits_{k=1}^\infty \bar c_k^2$ is bounded, the sum converges absolutely. Therefore, $S(t)$ converges uniformly on $[a,b]$. $\blacksquare$

Similarly, considering the estimate of $B_k$ and representation of $u_k(t)$ for $0\le t\le a$ we will have
\begin{eqnarray}\label{22}
    |u_k(t)|&\le& \left(c_{10}+\frac{c_{11}}{\lambda_k}\right)\frac{|\varphi_k|t}{1+\lambda_k t^\alpha}+\frac{c_{12}|\psi_k|t}{1+\lambda_kt^\alpha}+\frac{c_{13}t}{1+\lambda_kt^\alpha}\int\limits_0^a\frac{(a-z)^{\alpha-1}|f_k(z)|}{1+\lambda_k (a-z)^\alpha}dz\nonumber\\
    &+&\frac{c_{14}t}{1+\lambda_kt^\alpha}\int\limits_a^\xi\frac{(\xi-z)^{\beta-1}|f_k(z)|}{1+\lambda_k(\xi-z)^\beta}dz+c_{15}\int\limits_0^t\frac{(t-z)^{\alpha-1}|f_k(z)|}{1+\lambda_k(t-z)^\alpha}dz,\,\,0\le t\le a.
\end{eqnarray}
Here $c_j\,\,(j=\overline{10,15})$ are positive constants independent of $\lambda_k$.

Then considering \eqref{8}, we will get
\begin{eqnarray}\label{23}
    |u(t,x)|&\le& \sum\limits_{k=1}^\infty\left(c_{10}+\frac{c_{11}}{\lambda_k}\right)\frac{|\varphi_k|t}{1+\lambda_k t^\alpha}+\sum\limits_{k=1}^\infty\frac{c_{12}|\psi_k|t}{1+\lambda_kt^\alpha}+\sum\limits_{k=1}^\infty\frac{c_{13}t}{1+\lambda_kt^\alpha}\int\limits_0^a\frac{(a-z)^{\alpha-1}|f_k(z)|}{1+\lambda_k (a-z)^\alpha}dz\nonumber\\
    &+&\sum\limits_{k=1}^\infty\frac{c_{14}t}{1+\lambda_kt^\alpha}\int\limits_a^\xi\frac{(\xi-z)^{\beta-1}|f_k(z)|}{1+\lambda_k(\xi-z)^\beta}dz+\sum\limits_{k=1}^\infty c_{15}\int\limits_0^t\frac{(t-z)^{\alpha-1}|f_k(z)|}{1+\lambda_k(t-z)^\alpha}dz,\, (t,x)\in\Omega_1,
\end{eqnarray}
\begin{eqnarray}\label{24}
    |u_{xx}(t,x)|&\le& \sum\limits_{k=1}^\infty\left(c_{10}+\frac{c_{11}}{\lambda_k}\right)\frac{\lambda_k|\varphi_k|t}{1+\lambda_k t^\alpha}+\sum\limits_{k=1}^\infty\frac{c_{12}\lambda_k|\psi_k|t}{1+\lambda_kt^\alpha}+\sum\limits_{k=1}^\infty\frac{c_{13}\lambda_kt}{1+\lambda_kt^\alpha}\int\limits_0^a\frac{(a-z)^{\alpha-1}|f_k(z)|}{1+\lambda_k (a-z)^\alpha}dz\nonumber\\
    &+&\sum\limits_{k=1}^\infty\frac{c_{14}\lambda_kt}{1+\lambda_kt^\alpha}\int\limits_a^\xi\frac{(\xi-z)^{\beta-1}|f_k(z)|}{1+\lambda_k(\xi-z)^\beta}dz+\sum\limits_{k=1}^\infty c_{15}\lambda_k\int\limits_0^t\frac{(t-z)^{\alpha-1}|f_k(z)|}{1+\lambda_k(t-z)^\alpha}dz,\,(t,x)\in\Omega_1.
\end{eqnarray}

To prove the uniform convergence of infinite series in \eqref{24} we will prove the following statements:

Let $\varphi(x)$ and $\psi(x)$ be real-valued functions defined on $[0,1]$ with Fourier sine coefficients $\varphi_k$ and $\psi_k$. Assume $\lambda_k = (k\pi)^2$, $\alpha \in (1,2)$, and $c_{10}, c_{11}, c_{12}$ are positive constants. 

{\bf Lemma 3.} If $\varphi(x)$ and $\psi(x)$ satisfy the boundary conditions $\varphi(0)=\varphi(1)=0$ and $\psi(0)=\psi(1)=0$, and belong to the Sobolev space $H^2(0,1)$, then the following series converge uniformly on $t \in [0, a]$:
$$ \bar S_1(t) = \sum_{k=1}^\infty \left(c_{10}+\frac{c_{11}}{\lambda_k}\right) \frac{\lambda_k|\varphi_k|t}{1+\lambda_k t^\alpha}, 
$$
$$
\bar S_2(t) = \sum_{k=1}^\infty \frac{c_{12}\lambda_k|\psi_k|t}{1+\lambda_k t^\alpha}. 
$$
{\sc Proof:} By the Weierstrass M-test, we must find a uniform bound for the $t$-dependent multiplier for $t \in [0, a]$. Let $g(t) = \dfrac{\lambda_k t}{1+\lambda_k t^\alpha}$. We find its maximum by setting the derivative to zero:
$$ 
g'(t) = \frac{\lambda_k(1+\lambda_k t^\alpha) - \lambda_k t(\alpha \lambda_k t^{\alpha-1})}{(1+\lambda_k t^\alpha)^2} = \frac{\lambda_k + (1-\alpha)\lambda_k^2 t^\alpha}{(1+\lambda_k t^\alpha)^2} = 0. 
$$
Solving for $t$, the maximum occurs at $t_{max} = \left( \frac{1}{(\alpha-1)\lambda_k} \right)^{\frac{1}{\alpha}}$. Since $\alpha > 1$, as $k \to \infty$, $t_{max} \to 0$, ensuring it falls within the interval $[0,a]$. Substituting $t_{max}$ back into $g(t)$ yields the peak value:
$$
g(t_{max}) = \frac{\lambda_k \left[ (\alpha-1)\lambda_k \right]^{-\frac{1}{\alpha}}}{1 + \frac{1}{\alpha-1}} = \frac{(\alpha-1)^{1-\frac{1}{\alpha}}}{\alpha} \lambda_k^{1-\frac{1}{\alpha}}. 
$$
Because $\lambda_k = k^2\pi^2$, the maximum value grows asymptotically as $O(k^{2-\frac{2}{\alpha}})$. Looking at $\bar S_1(t)$ and $\bar S_2(t)$, the term $\left(c_{10}+\dfrac{c_{11}}{\lambda_k}\right)$ is bounded by a constant for large $k$. Therefore, for uniform convergence, we require:
$$
\sum_{k=1}^\infty k^{2-\frac{2}{\alpha}} |\varphi_k| < \infty \quad \text{and} \quad \sum_{k=1}^\infty k^{2-\frac{2}{\alpha}} |\psi_k| < \infty .
$$
Given that $\varphi(0)=\varphi(1)=0$, applying integration by parts twice to the Fourier definition yields $\varphi_k = -\dfrac{s_k}{k^2\pi^2}$, where $s_k$ are the Fourier sine coefficients of the second derivative $\varphi''(x)$. If $\varphi \in H^2(0,1)$, then $\sum s_k^2 < \infty$. Substituting this into our requirement and applying the Cauchy-Schwarz inequality:
$$
\sum_{k=1}^\infty k^{2-\frac{2}{\alpha}} \frac{|s_k|}{k^2\pi^2} = \frac{1}{\pi^2} \sum_{k=1}^\infty k^{-\frac{2}{\alpha}} |s_k| \le \frac{1}{\pi^2} \left( \sum_{k=1}^\infty s_k^2 \right)^{\frac{1}{2}} \left( \sum_{k=1}^\infty k^{-\frac{4}{\alpha}} \right)^{\frac{1}{2}}. 
$$
Since $1 < \alpha < 2$, we have $\frac{4}{\alpha} > 2$. Thus, the series $\sum k^{-4/\alpha}$ converges by the p-series test, guaranteeing the absolute and uniform convergence of $\bar S_1(t)$ and $\bar S_2(t)$. $\blacksquare$

Let $f_k(t)$ denote the spatial Fourier sine coefficients of $f(t,x)$ for $x \in [0,1]$. Let $\lambda_k = (k\pi)^2$, $\beta \in (0,1)$, $\alpha \in (1,2)$, and $\xi \in (a,b)$.

{\bf Lemma 4.}  If $f(t,x)$ satisfies the boundary conditions $f(t,0) = f(t,1) = 0$ for all $t$, and $f(t, \cdot) \in H^1(0,1)$, then three series then the following series converge uniformly on $t \in [0, a]$:
$$
\bar S_3= \sum\limits_{k=1}^\infty\frac{c_{13}\lambda_kt}{1+\lambda_kt^\alpha}\int\limits_0^a\frac{(a-z)^{\alpha-1}|f_k(z)|}{1+\lambda_k (a-z)^\alpha}dz,
$$
$$
\bar S_4=\sum\limits_{k=1}^\infty\frac{c_{14}\lambda_kt}{1+\lambda_kt^\alpha}\int\limits_a^\xi\frac{(\xi-z)^{\beta-1}|f_k(z)|}{1+\lambda_k(\xi-z)^\beta}dz,
$$
$$
\bar S_5= \sum\limits_{k=1}^\infty c_{15}\lambda_k\int\limits_0^t\frac{(t-z)^{\alpha-1}|f_k(z)|}{1+\lambda_k(t-z)^\alpha}dz.
$$

{\sc Proof:} We must check the required convergence rate of $M_k = \sup\limits_t |f_k(t)|$ for the three remaining series. For $\bar S_3(t)$ and $\bar S_4(t)$: These series feature the same external multiplier $g(t) = \dfrac{\lambda_k t}{1+\lambda_k t^\alpha}$, which we established is bounded by $O(k^{2-\frac{2}{\alpha}})$. The integral portion of $\bar S_3$ evaluates to:
$$
\int\limits_0^a \frac{(a-z)^{\alpha-1}|f_k(z)|}{1+\lambda_k(a-z)^\alpha} dz \le M_k \frac{1}{\alpha \lambda_k} \ln(1+\lambda_k a^\alpha) = O\left(M_k \frac{\ln k}{k^2}\right). 
$$
Combining the multiplier and the integral, the overall terms of $\bar S_3$ (and similarly $\bar S_4$) are bounded by:
$$
O\left(k^{2-\frac{2}{\alpha}}\right) \cdot O\left(M_k \frac{\ln k}{k^2}\right) = O\left( M_k k^{-\frac{2}{\alpha}} \ln k \right). 
$$
Because $\alpha < 2$, the exponent $-\frac{2}{\alpha} < -1$. Thus, $\bar S_3$ and $\bar S_4$ will converge absolutely as long as $M_k$ is merely bounded (i.e., $f(t,x)$ is continuous). For $\bar S_5(t)$: This series does not contain the external $t$ multiplier. Its terms are bounded solely by the evaluation of its integral:
$$
\left| c_{15}\lambda_k\int\limits_0^t\frac{(t-z)^{\alpha-1}|f_k(z)|}{1+\lambda_k(t-z)^\alpha}dz \right| \le c_{15} \lambda_k M_k \frac{1}{\alpha \lambda_k} \ln(1+\lambda_k t^\alpha) = O(M_k \ln k). 
$$
This imposes the strictest requirement on $f(t,x)$: we must have $\sum\limits_{k=1}^\infty M_k \ln k < \infty$. Because $f(t,0)=f(t,1)=0$, applying integration by parts once yields $M_k = \dfrac{\tilde C_k}{k\pi}$, where $\tilde C_k$ are the uniformly bounded Fourier cosine coefficients of $\dfrac{\partial f}{\partial x}$. Because $\dfrac{\partial f}{\partial x} \in L^2$, we know $\sum \tilde C_k^2 < \infty$. Testing the required sum via Cauchy-Schwarz:
$$
\sum_{k=1}^\infty M_k \ln k \le \frac{1}{\pi} \sum_{k=1}^\infty \tilde C_k \frac{\ln k}{k} \le \frac{1}{\pi} \left( \sum_{k=1}^\infty \tilde C_k^2 \right)^{\frac{1}{2}} \left( \sum_{k=1}^\infty \frac{\ln^2 k}{k^2} \right)^{\frac{1}{2}}. 
$$
Since $\sum \dfrac{\ln^2 k}{k^2}$ converges, the series $\bar S_5(t)$ converges uniformly, which trivially guarantees the convergence of $\bar S_3(t)$ and $\bar S_4(t)$ as well. $\blacksquare$

Using \eqref{17}, we evaluate the function $h(x)$
      		\begin{equation*}
			\begin{aligned}
				&\left| h\left(x \right) \right|\le \sum\limits_{k=1}^\infty f_k(a)+\sum\limits_{k=1}^\infty \left(\lambda_k a^{\alpha-1}+\frac{c_4}{\lambda_k}\right)\frac{c_{16}|\varphi_k|}{1+\lambda_k a^\alpha}+\sum\limits_{k=1}^\infty \left(c_1\lambda_k +\frac{c_3c_{16}}{1+\lambda_ka^\alpha}\right)|\psi_k|\\
                &+\sum\limits_{k=1}^\infty c_{16}\int\limits_0^a \frac{(a-z)^{\alpha-1}}{1+\lambda_k(a-z)^\alpha}\left(\frac{c_5(a-z)}{1+\lambda_k a^\alpha}+1\right)|f_k(z)|dz+
                \sum\limits_{k=1}^\infty \left(c_2\lambda_k+\frac{c_6c_{16}}{1+\lambda_ka^\alpha}\right)\int\limits_a^\xi \frac{(\xi-z)^{\beta-1}}{1+\lambda_k(\xi-z)^\beta}|f_k(z)|dz,
			\end{aligned}
		\end{equation*}
		where $c_{16}$ is a positive constant independent of $\lambda_k$.

    \textbf{Theorem 1.}
Let $1<\alpha\le 4/3$. If
\begin{itemize}
    \item $f(\cdot,x)$ is continuous  for all $t \in [0,b]$ such that $f(t,0) = f(t,1) = 0$, and $f(t, \cdot)$ is absolutely continuous with its partial derivative $\frac{\partial f}{\partial x} \in L^2(0,1)$;
    \item $\psi(x)\in AC^2[0,1],\,\psi''' \in L^2(0,1)$ such that $\psi(0) = \psi(1) = \psi''(0) = \psi''(1) = 0$;
    \item $\varphi(x)\in AC^1[0,1],\,\varphi'' \in L^2(0,1)$ such that $\varphi(0) = \varphi(1) = 0$,
\end{itemize}
then there exists a unique solution of problem 1, represented by $u(t,x)$ \eqref{8} and \eqref{h_k}.

{\bf Remark 1.} If $4/3<\alpha < 2$, the statement of Theorem 1 will be valid with the condition $E_{\alpha, 2}\left[-\lambda_ka^{\alpha}\right]\neq 0$ for all $k\in\mathbb{N}$.

\section{Inverse problem on finding unknown part of position/concentration transmitting condition}

\textbf{Problem 2.} We aim to find functions $u(t,x) \in C(\bar{\Omega})\cap C^2(\Omega_1\cup\Omega_2)$, $u(\cdot,x) \in AC^2[0,a)\cup AC^1[a,b),$ $\bar h(x)\in C[0,1]$, which satisfies the equation \eqref{1}, together with the initial condition \eqref{2}, boundary conditions \eqref{3}, the over-determination condition \eqref{4}
and the transmitting conditions:
\begin{equation} \label{ulash1}
    \lim_{t \rightarrow a+0} { }^C D_{at}^{\beta} u(t, x)=u_t(a-0,x), \quad 0< x< 1,
\end{equation}
\begin{equation} \label{ulash2}
    u(a+0,x)=u(a-0,x)+\bar h(x), \quad 0\le x\le 1.
\end{equation}

This problem is studied in the same way as Problem 1 considered above. The function $u(x,t)$, defined by formulas \eqref{13} and \eqref{14}, must satisfy the transmitting conditions \eqref{ulash1} and \eqref{ulash2}. In this case, from the first transmitting condition \eqref{ulash1}, \eqref{13} and \eqref{14}, assuming that  $E_{\alpha,1}[\lambda_k a^\alpha]\neq 0$ for all $k\in\mathbb{N}$, we will find $B_k$:
\begin{equation} \label{B_k}
    \begin{gathered}
    B_k=\frac{f_k(a)+\lambda_k\left(\varphi_k a^{\alpha-1}E_{\alpha,\alpha}\left[-\lambda_k a^{\alpha}\right]- C_k\right) -\int\limits_0^a (a-z)^{\alpha-2} E_{\alpha, \alpha-1}\left[-\lambda_k (a-z)^{\alpha}\right] f_k(z) dz}{E_{\alpha,1}[\lambda_k a^\alpha]}.
\end{gathered}
\end{equation}

Using the second transmitting condition \eqref{ulash2}, we will find $h_k$:
\begin{equation} \label{h_k2}
    \begin{gathered}
\bar h_k=C_k-\varphi_k E_{\alpha,1}\left[-\lambda_k a^{\alpha}\right]-B_k a E_{\alpha,2}\left[-\lambda_k a^{\alpha}\right]
-\int_0^a (a-z)^{\alpha-1} E_{\alpha, \alpha}\left[-\lambda_k (a-z)^{\alpha}\right] f_k(z) dz.
\end{gathered}
\end{equation}

With $B_k$ and $\bar{h}_k$ given by \eqref{B_k}-\eqref{h_k2}, the candidate
solution is again
\[
  u(t,x) = \sum_{k=1}^{\infty} u_k(t)\sin(k\pi x),
  \qquad
  \bar{h}(x) = \sum_{k=1}^{\infty} \bar{h}_k \sin(k\pi x),
\]
where $u_k(t)$ is the same piecewise formula \eqref{13}, \eqref{14} as in Problem~1, but
with the new $B_k$ coefficient \eqref{B_k}.
 
Using the Mittag-Leffler estimate $|E_{a,b}(z)| \le c/(1+|z|)$ for $z \le 0$,
$0 < a < 2$ (see \eqref{MLE}), one derives from \eqref{B_k} the bound
 
\begin{eqnarray}\label{eq:Bk_P2_bound}
  |B_k| \;&\le&\;
  \frac{c_1}{|E_{\alpha,1}(-\lambda_k a^\alpha)|}\nonumber\\
  &\times&\left[
    \frac{|f_k(a)|}{\lambda_k}
    + \frac{c_2|\psi_k|}{1+\lambda_k(\xi-a)^\beta}
    + c_3\int_a^{\xi}
        \frac{(\xi-z)^{\beta-1}|f_k(z)|}{1+\lambda_k(\xi-z)^\beta}\,dz
    + c_4\int_0^a
        \frac{(a-z)^{\alpha-2}|f_k(z)|}{1+\lambda_k(a-z)^\alpha}\,dz
  \right],
\end{eqnarray}
provided the denominators $E_{\alpha,1}(-\lambda_k a^\alpha)$ are bounded away from zero
uniformly in $k$.
 
Similarly, from \eqref{h_k2} one obtains
\begin{equation}\label{eq:hbk_P2_bound}
  |\bar{h}_k| \;\le\;
  c_5 |\psi_k| + \frac{c_6}{\lambda_k}|\varphi_k|
  + c_7\int_0^a \frac{(a-z)^{\alpha-1}|f_k(z)|}{1+\lambda_k(a-z)^\alpha}\,dz
  + c_8\int_a^{\xi}\frac{(\xi-z)^{\beta-1}|f_k(z)|}{1+\lambda_k(\xi-z)^\beta}\,dz,
\end{equation}
with constants $c_j > 0$ independent of $k$ and $\lambda_k$.
 
The uniform convergence of the series for $u$, $u_{xx}$, ${}^C\!D_{0t}^\alpha u$,
${}^C\!D_{at}^\beta u$, and $\bar{h}$ follows from the same four lemmas used in
Problem~1 (Lemmas~1-4), once \eqref{eq:Bk_P2_bound}--\eqref{eq:hbk_P2_bound} are
substituted for the corresponding bounds in Problem~1.
We therefore state the main result directly.

{\bf Theorem 2.}  Let $1 < \alpha \le \tfrac{4}{3}$, $0 < \beta < 1$, $0 < a < b$, and $\xi \in (a,b)$.
  Suppose that the following conditions hold:
  \begin{enumerate}
    \item\label{it:f_P2}
      $f(t,x)$ is continuous on $[0,b]\times[0,1]$ with $f(t,0)=f(t,1)=0$ for all
      $t\in[0,b]$, and $f(t,\cdot)\in H^1_0(0,1)$ is absolutely continuous with
      $\partial_x f(t,\cdot)\in L^2(0,1)$;
    \item\label{it:psi_P2}
      $\psi(x)\in AC^2[0,1]$, $\psi'''\in L^2(0,1)$, and
      $\psi(0)=\psi(1)=\psi''(0)=\psi''(1)=0$;
    \item\label{it:phi_P2}
      $\varphi(x)\in AC^1[0,1]$, $\varphi''\in L^2(0,1)$, and $\varphi(0)=\varphi(1)=0$;
    \item\label{it:nonzero_P2}
      $E_{\alpha,1}(-\lambda_k a^\alpha)\neq 0$ for all $k\in\mathbb{N}$.
  \end{enumerate}
  Then Problem~\emph{2} has a unique solution
  \[
    \bigl(u(t,x),\;\bar{h}(x)\bigr)
    \;\in\;
    \Bigl[C(\bar\Omega)\cap C^2(\Omega_1\cup\Omega_2),\;
          u(\cdot,x)\in AC^2[0,a)\cup AC^1[a,b)\Bigr]
    \;\times\; C[0,1],
  \]
  represented by the uniformly convergent series
  \begin{align*}
    u(t,x) = \sum_{k=1}^{\infty} u_k(t)\sin(k\pi x), \,\,\,\,
    \bar{h}(x) = \sum_{k=1}^{\infty} \bar{h}_k\,\sin(k\pi x).
  \end{align*}
 
{\sc Proof:}
  The argument follows the same four-step architecture as the proof of Theorem~1, with
the two transmitting conditions interchanging their roles.
We highlight only the points where Problem~2 differs.
 
\smallskip
\noindent\textbf{Step 1. Reduction to decoupled ODEs.}
Substituting the Fourier ansatz $u(t,x)=\sum\limits_k u_k(t)\sin(k\pi x)$ into the governing
equation and the boundary conditions yields the same pair of fractional Cauchy problems on $(0,a)$ and $(a,b)$ as in Problem~1,
with the same solution templates \eqref{13}, \eqref{14}.
 
\smallskip
\noindent\textbf{Step 2. Determination of $C_k$.}
The overdetermination condition $u(\xi,x)=\psi(x)$ acts on the diffusion piece and gives $C_k$ by formula, identical to
Problem~1.
 
\smallskip
\noindent\textbf{Step 3. Determination of $B_k$ from the velocity transmitting
condition~\eqref{ulash1}.}
This is where Problem~2 departs from Problem~1.
Evaluating $\lim_{t\to a^+}{}^C\!D_{at}^\beta u_k(t)$ from the diffusion piece and $u_k'(a^-)$ from the wave piece, the
condition $\lim\limits_{t\to a^+}{}^C\!D_{at}^\beta u_k = u_k'(a^-)$ gives
\[
  f_k(a) - \lambda_k C_k
  \;=\;
  -\lambda_k\varphi_k a^{\alpha-1}E_{\alpha,\alpha}(-\lambda_k a^\alpha)
  + B_k E_{\alpha,1}(-\lambda_k a^\alpha)
  + \int_0^a (a-z)^{\alpha-2}E_{\alpha,\alpha-1}(-\lambda_k(a-z)^\alpha)f_k(z)\,dz.
\]
Solving for $B_k$ (using hypothesis~\eqref{it:nonzero_P2}) yields \eqref{B_k}.
Observe that in Problem~1 it was the function $E_{\alpha,2}$ that appeared in the
denominator (from the position continuity condition); here it is $E_{\alpha,1}$
(from the velocity condition).
 
\smallskip
\noindent\textbf{Step 4. Determination of $\bar{h}_k$ from the position transmitting
condition~\eqref{ulash2}.}
The condition $u_k(a^+) = u_k(a^-) + \bar{h}_k$ gives immediately, from
\eqref{13}-\eqref{14} with $E_{\beta,1}(0)=1$,
\[
  \bar{h}_k
  = C_k - \varphi_k E_{\alpha,1}(-\lambda_k a^\alpha)
    - B_k\,a\,E_{\alpha,2}(-\lambda_k a^\alpha)
    - \int_0^a (a-z)^{\alpha-1}E_{\alpha,\alpha}(-\lambda_k(a-z)^\alpha)f_k(z)\,dz,
\]
which is \eqref{h_k2}.
Note the structural duality: in Problem~1, $h_k$ was recovered from the velocity
condition (Step~3 of that proof) once $B_k$ was found from the position condition;
here the roles of the two transmitting conditions are exchanged.
 
\smallskip
\noindent\textbf{Step 5. Uniform convergence.}
We must verify that the series
\[
  \sum_{k=1}^\infty u_k(t)\sin(k\pi x),\quad
  \sum_{k=1}^\infty u_{k,xx}(t)\sin(k\pi x),\quad
  \sum_{k=1}^\infty \bar{h}_k\sin(k\pi x)
\]
converge uniformly under hypotheses \eqref{it:f_P2}--\eqref{it:phi_P2}.
 
The bounds on $C_k$ are identical to Problem~1.
The bounds on $B_k$ from \eqref{eq:Bk_P2_bound} contain an extra term involving
$f_k(a)/(\lambda_k)$ and $(a-z)^{\alpha-2}$ instead of $(a-z)^{\alpha-1}$.
Because $\alpha>1$ one has $\alpha-2>-1$, so $(a-z)^{\alpha-2}$ is integrable on
$[0,a]$, and the corresponding integral is $O(M_k/\lambda_k)$, which is summable
under hypothesis~\eqref{it:f_P2} by the same Cauchy–Schwarz argument used in
Lemma~1.
The Weierstrass $M$-test then applies exactly as in the proof of Lemmas~1-2 to
give uniform convergence of the series on $\bar\Omega_2$.
 
The estimate for $u_k(t)$ on $[0,a]$ takes the same form as \eqref{14}
but with the new $B_k$.
The factor $\lambda_k t/(1+\lambda_k t^\alpha)$ is bounded by $O(k^{2-2/\alpha})$
(Lemma~3), and because $4/\alpha > 2$ for $\alpha < 2$, the resulting $p$-series
converges.
The integral contributions involving $f_k$ are handled by Lemma~4, exactly as before.
 
From \eqref{eq:hbk_P2_bound}, the dominant term is $c_5|\psi_k|$, which requires
$\sum\limits_k k^2|\psi_k|<\infty$; this follows from $\psi\in H^3(0,1)$ exactly as in
Lemma~2.
The term $c_6|\varphi_k|/\lambda_k$ requires $\sum\limits_k|\varphi_k|<\infty$, which holds
if $\varphi\in H^1(0,1)$, guaranteed by hypothesis~\eqref{it:phi_P2}.
The integral terms are treated by Lemma~1 and Lemma~4.
 
\smallskip
\noindent\textbf{Step 6. Uniqueness.}
Suppose $(u^{(1)}, \bar{h}^{(1)})$ and $(u^{(2)}, \bar{h}^{(2)})$ are two solutions.
Setting $w = u^{(1)}-u^{(2)}$ and $g = \bar{h}^{(1)}-\bar{h}^{(2)}$, one finds that
$w$ satisfies the homogeneous version of Problem~2 (zero forcing, initial, boundary,
and overdetermination data), together with a jump $g(x)$ in the position transmitting
condition.
Projecting onto $\sin(k\pi x)$ gives $w_k(\xi)=0$ for each $k$.
In the diffusion region $(a,b)$, the equation ${}^C\!D_{at}^\beta w_k + \lambda_k w_k = 0$
with $w_k(\xi)=0$ admits only the trivial solution $w_k\equiv 0$ on $[a,b]$, since
$E_{\beta,1}(-\lambda_k(\xi-a)^\beta)\neq 0$ for all $0<\beta<1$.
Hence $C_k = 0$ for all $k$.
In the wave region $(0,a)$, the homogeneous system then gives $B_k=0$ (by
hypothesis~\eqref{it:nonzero_P2}) and $w_k\equiv 0$ on $[0,a]$.
Finally, $g_k = 0$ follows from \eqref{h_k2}.
Thus $w\equiv 0$ and $g\equiv 0$, establishing uniqueness.

{\bf Remark 2.} Condition~\eqref{it:nonzero_P2} replaces the condition $E_{\alpha,2}(-\lambda_k
  a^\alpha)\neq 0$ that appeared in Problem~1 (Remark~1), and it plays the same
  structural role: it ensures that the formula for $B_k$ in \eqref{B_k} is
  non-singular.
  Note that $E_{\alpha,1}(z)$ is the Mittag-Leffler function with second parameter
  $b=1$; for $1<\alpha\le 4/3$ and $z<0$, it is known to have at most finitely many
  real zeros~\cite{PS2011}, so condition~\eqref{it:nonzero_P2} fails for at
  most finitely many $k$.
  If it fails for some $k=k_0$, that mode is excluded from the spectral expansion, and
  the result holds for the reduced system.

\section{An explicit example and numerical illustration}

\subsection{Choice of data}

To illustrate the solution formulas derived in Problem~1 we choose data that make the
Fourier series collapse to a single mode, yielding a fully closed-form expression.
Set
\begin{equation}\label{eq:data}
  \varphi(x) = \sin(\pi x), \quad
  \psi(x)    = \sin(\pi x), \quad
  f(t,x)     = 0,
\end{equation}
and fix the domain parameters
\[
  a = 0.5, \quad b = 1, \quad \xi = 0.75.
\]

Because $\{\sin(k\pi x)\}_{k=1}^{\infty}$ is an orthonormal basis of $L^2(0,1)$, the Fourier
coefficients of the data are
\[
  \varphi_k = \psi_k = \delta_{k,1}, \qquad f_k(t) \equiv 0 \quad \text{for all } k \geq 1.
\]
Only the $k=1$ eigenvalue $\lambda_1 = \pi^2$ contributes, and the solution takes the
finite form
\begin{equation}\label{eq:sol_finite}
  u(t,x) = u_1(t)\,\sin(\pi x), \qquad h(x) = h_1\,\sin(\pi x),
\end{equation}
where $u_1(t)$ and $h_1$ are computed below.

\subsection{Explicit computation of the Fourier coefficients}

With $f \equiv 0$, all integral terms in formulas~(14)--(18) of the main text vanish.
Setting $\lambda = \lambda_1 = \pi^2$ to lighten notation, the three key constants are
obtained in closed form.

\paragraph{Step 1. Determine $C_1$ from the overdetermination condition.}
Formula~(16) with $f_k \equiv 0$ gives
\begin{equation}\label{eq:C1}
  C_1 = \frac{\psi_1}{E_{\beta,1}\!\left[-\lambda({\xi-a})^{\beta}\right]}.
\end{equation}

\paragraph{Step 2. Determine $B_1$ from the position transmitting condition.}
Formula~(17) with $f_k \equiv 0$ gives
\begin{equation}\label{eq:B1}
  B_1 = \frac{C_1 - \varphi_1 E_{\alpha,1}(-\lambda a^{\alpha})}
             {a\, E_{\alpha,2}(-\lambda a^{\alpha})}.
\end{equation}

\paragraph{Step 3. Determine $h_1$ from the velocity transmitting condition.}
Formula~(18) with $f_k \equiv 0$ gives
\begin{equation}\label{eq:h1}
  h_1 = -\lambda C_1 - B_1 E_{\alpha,1}(-\lambda a^{\alpha})
        + \lambda\,\varphi_1\, a^{\alpha-1}\, E_{\alpha,\alpha}(-\lambda a^{\alpha}).
\end{equation}

\paragraph{Piecewise solution.}
\begin{equation}\label{eq:u1_pw}
  u_1(t) =
  \begin{cases}
    \varphi_1\, E_{\alpha,1}(-\lambda t^{\alpha})
    + B_1\, t\, E_{\alpha,2}(-\lambda t^{\alpha}), & 0 \leq t \leq a, \\[6pt]
    C_1\, E_{\beta,1}\!\left[-\lambda(t-a)^{\beta}\right], & a < t \leq b.
  \end{cases}
\end{equation}

{\bf Remark 3.}
The condition $E_{\alpha,2}(-\lambda a^{\alpha}) \neq 0$ is required for $B_1$ to be
well-defined; by Theorem~1, this is guaranteed when $1 < \alpha \leq \frac{4}{3}$.
For the values used below ($\alpha \in \{1.2, 1.3, 1.5, 1.8\}$), we verify numerically
that $E_{\alpha,2}(-\lambda a^{\alpha}) > 0$ in all cases.

\subsection{Numerical values for selected fractional orders}

Table~\ref{tab:numerics} lists the intermediate and final quantities for four
representative pairs $(\alpha,\beta)$, computed via the algorithm of
Garrappa~\cite{Garrappa2015} for the two-parameter Mittag-Leffler function.

\begin{table}[h]
\centering
\caption{Numerical values of the key Mittag-Leffler evaluations and the resulting
  constants $C_1$, $B_1$, $h_1$ for $\lambda = \pi^2$, $a = 0.5$, $\xi = 0.75$.}
\label{tab:numerics}
\renewcommand{\arraystretch}{1.25}
\begin{tabular}{cccccccc}
\hline
$\alpha$ & $\beta$ &
$E_{\alpha,1}(-\lambda a^\alpha)$ &
$E_{\alpha,2}(-\lambda a^\alpha)$ &
$E_{\beta,1}(-\lambda(\xi-a)^\beta)$ &
$C_1$ & $B_1$ &$h_1$ \\ \hline
$1.8$ & $0.3$ & $-0.17677$ & $0.51294$ & $0.10813$ & $9.2479$  & $36.748$  & $-82.106$ \\
$1.5$ & $0.5$ & $-0.23376$ & $0.33531$ & $0.11211$ & $8.9196$  & $54.596$  & $-74.264$ \\
$1.2$ & $0.7$ & $-0.07366$ & $0.23358$ & $0.10763$ & $9.2909$  & $80.182$ & $-85.734$ \\
$1.3$ & $0.4$ & $-0.13134$ & $0.25789$ & $0.11099$ & $9.0097$  & $70.892$  & $-79.425$ \\
\hline
\end{tabular}

\end{table}

\subsection{Verification of conditions}

\begin{enumerate}
  \item \textbf{Overdetermination:}
    $u_1(\xi) = C_1\,E_{\beta,1}(-\lambda(\xi-a)^\beta)
              = C_1\cdot E_{\beta,1}(-\lambda(\xi-a)^\beta) = \psi_1 = 1$
    holds exactly by construction~\eqref{eq:C1}.

  \item \textbf{Position continuity at $t=a$:}
    \[
      u_1(a^-)
        = \varphi_1 E_{\alpha,1}(-\lambda a^\alpha) + B_1\,a\,E_{\alpha,2}(-\lambda a^\alpha)
        = C_1\,E_{\beta,1}(0) = C_1
        = u_1(a^+),
    \]
    which follows immediately from~\eqref{eq:B1} and $E_{\beta,1}(0)=1$.

  \item \textbf{Initial condition:}
    $u_1(0) = \varphi_1\,E_{\alpha,1}(0) = \varphi_1 = 1$, since $E_{\alpha,1}(0)=1$.
\end{enumerate}

\subsection{Qualitative behaviour}

\begin{figure}[h]
  \centering
  \includegraphics[width=\textwidth]{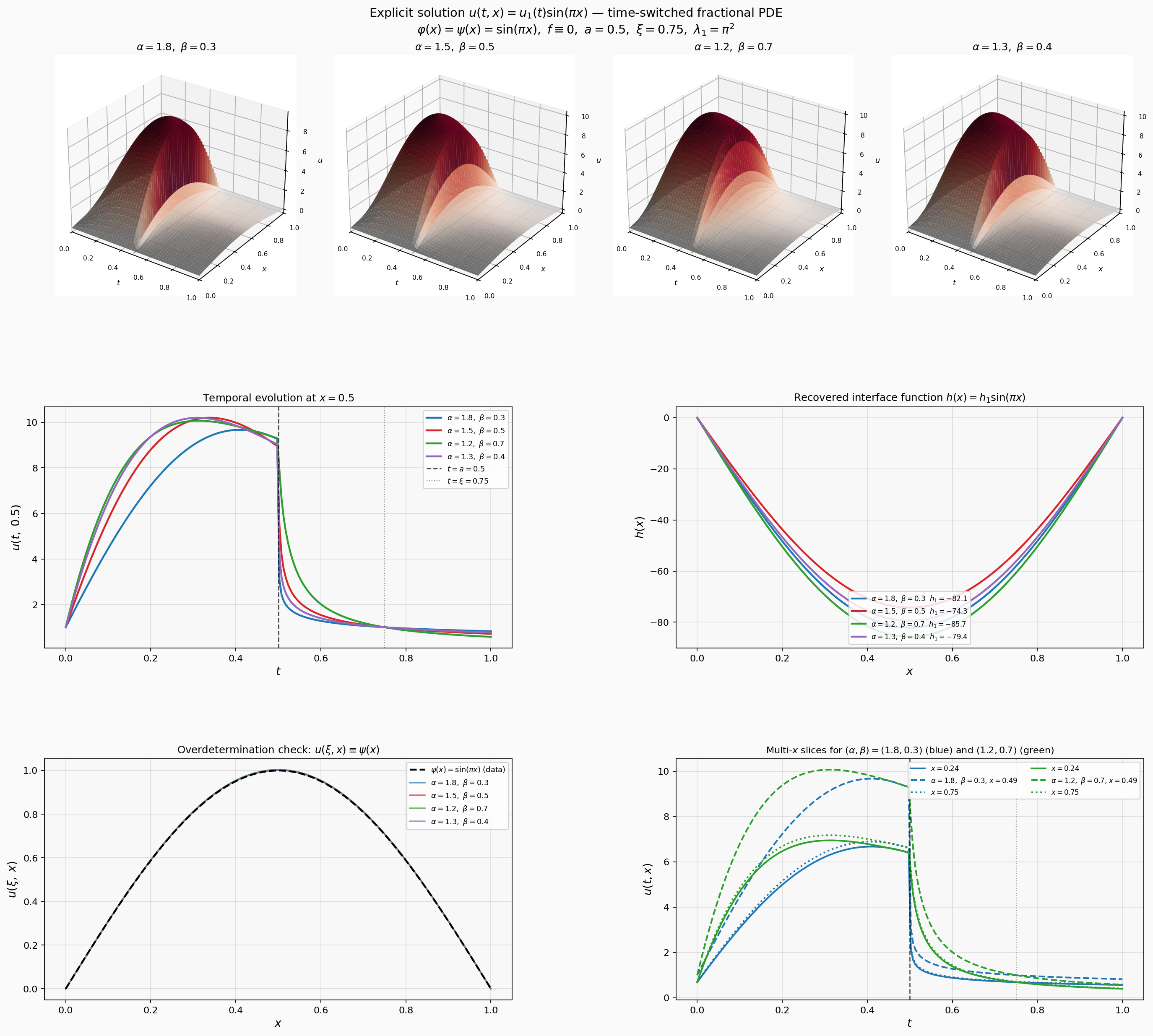}
  \caption{%
    Explicit solution $u(t,x)=u_1(t)\sin(\pi x)$ of Problem~1 for four pairs
    $(\alpha,\beta)$.\\\textbf{Top row:} 3-D solution surfaces; the dashed vertical plane marks $t=a=0.5$.\\
    \textbf{Middle left:} temporal traces $u(t,0.5)$ for all four cases; the vertical
    dashed line marks the switching time $a=0.5$ and the dotted line marks the
    overdetermination instant $\xi=0.75$.\\
    \textbf{Middle right:} recovered unknown $h(x)=h_1\sin(\pi x)$.\\
    \textbf{Lower left:} spatial profiles $u(\xi,x)$ overlaid on the exact data
    $\psi(x)=\sin(\pi x)$, confirming the overdetermination condition.\\
    \textbf{Lower right:} multi-$x$ slices for the cases $(\alpha,\beta)=(1.8,0.3)$
    (blue) and $(1.2,0.7)$ (green), illustrating the dependence of the solution
    amplitude on the spatial location.%
  }
  \label{fig:solution}
\end{figure}

Figure~\ref{fig:solution} displays the solution surface $u(t,x)$ and several derived
quantities for the four $(\alpha,\beta)$ pairs in Table~\ref{tab:numerics}.
The following features are visible.

\begin{enumerate}
  \item \textbf{Wave-phase oscillation ($0 \leq t \leq a$).}
    For larger $\alpha$ (closer to~$2$) the Mittag-Leffler function
    $E_{\alpha,1}(-\lambda t^\alpha)$ exhibits more pronounced oscillatory decay, while
    for $\alpha$ close to~$1$ the behaviour is closer to the exponential, reflecting the
    transition to the diffusion regime.

  \item \textbf{Diffusion-phase monotone decay ($a < t \leq b$).}
    After switching, $u_1(t) = C_1 E_{\beta,1}(-\lambda(t-a)^\beta)$ decays monotonically
    because $E_{\beta,1}(-\lambda\tau^\beta)$ is completely monotone for $0<\beta<1$.
    Smaller $\beta$ produces heavier-tailed (slower) subdiffusive decay.

  \item \textbf{Continuity and kink at $t=a$.}
    The solution is continuous at $t=a$ by design, but the velocity
    $\partial_t u$ has a jump proportional to $h_1\sin(\pi x)$, which is the
    quantity recovered by the inverse problem.

  \item \textbf{Recovered interface function $h(x)$.}
    The function $h(x) = h_1\sin(\pi x)$ is negative for all tested parameter pairs,
    indicating that the diffusion-side fractional flux is smaller than the wave-side
    velocity at $t=a$. The magnitude $|h_1|$ is largest when both $\alpha$ is close
    to~$1$ and $\beta$ is close to~$1$ (case $\alpha=1.2,\beta=0.7$, $|h_1|\approx 85.7$),
    consistent with the fact that steeper fractional gradients require a larger
    correction at the interface.

  \item \textbf{Overdetermination check.}
    The spatial profile $u(\xi,x)$ coincides with $\psi(x)=\sin(\pi x)$ for all four
    parameter sets (Figure~\ref{fig:solution}, lower-left panel), confirming the
    algebraic consistency of the explicit formulae.
\end{enumerate}

\end{document}